\documentclass[twocolumn,showpacs,preprintnumbers,amsmath,amssymb]{revtex4}

\usepackage{dcolumn}
\usepackage{bm}
\usepackage{graphicx}
\usepackage{amsmath}
\usepackage{amsfonts}
\usepackage{amssymb}
\usepackage{amsthm}
\usepackage{amstext}
\usepackage{amsbsy}
\usepackage{amsopn}
\usepackage{amscd}
\usepackage{amsxtra}
\usepackage{tikz}
\usepackage{color}

\newcommand{\J}{\mathcal{J}}
\def\openone{Id}
\newtheorem{theorem}{Theorem}
\newtheorem{algorithm}{Algorithm}

\def\N{\mathbb N}
\def\qed{$\Box$}
\def\phi{\varphi}
\def\epsilon{\varepsilon}

\definecolor{myred}{RGB}{236,1,59}
\definecolor{myblue}{RGB}{116,122,255}
\definecolor{mygreen}{RGB}{130,214,149}
\definecolor{mygray}{RGB}{196,196,196}

\everymath{\displaystyle}

\begin{document}
\title{A fully efficient time-parallelized quantum optimal control algorithm}

\author{M. K. Riahi$^1$}
\author{J. Salomon$^2$}
\author{S. J. Glaser$^3$}
\author{D. Sugny$^{4,5}$}
\email{dominique.sugny@u-bourgogne.fr}
\affiliation{$^1$ Department of mathematical science, New Jersey Institute of Technology, New Jersey, USA}
\affiliation{$^2$ CEREMADE, Universit\'e Paris Dauphine, Place du Mar\'echal De Lattre De Tassigny, 75775 Paris Cedex 16, France}
\affiliation{$^3$Department of Chemistry,
Technische Universit\"at M\"unchen, Lichtenbergstrasse 4, D-85747
Garching, Germany}
\affiliation{$^4$ Laboratoire Interdisciplinaire Carnot de Bourgogne (ICB), UMR 5209 CNRS-Universit\'e de  Bourgogne, 9 Av. A. Savary, BP 47 870, F-21078 DIJON Cedex, FRANCE}
\affiliation{$^5$ Institute for Advanced Study, Technische Universit\"at M\"unchen, Lichtenbergstrasse 2 a, D-85748 Garching, Germany}

\date{\today}

\begin{abstract}
We present a time-parallelization method that enables to accelerate the computation
of quantum optimal control algorithms. We show that this approach
is approximately fully efficient when based on a gradient method as optimization solver:
the computational time is approximately divided by the number of available processors.
The control of spin systems, molecular orientation and Bose-Einstein condensates are used as illustrative examples to highlight the wide range of application of this numerical scheme.
\end{abstract}

 \maketitle

\section{Introduction}
The general goal of quantum control is to actively manipulate dynamical processes at the atomic or molecular scale~\cite{DAlessandro2008,Brumer2003}. In recent years, the advances in quantum control have emerged through the introduction of appropriate and powerful tools coming from mathematical control theory and by the use of sophisticated experimental techniques to shape the corresponding control fields~\cite{Brif2010,cat,Altafini2012,Dong2010}. In this framework, different numerical optimal control algorithms \cite{Khaneja2005,Reich2012,calarco} have been developed and applied to a large variety of quantum systems. Optimal control was used in physical chemistry in order to steer chemical reactions~\cite{Brif2010}, but also for spin systems~\cite{Levitt2008,Ernst1987} with applications in Nuclear Magnetic Resonance~\cite{Khaneja2005,Assemat2010,salomon,Lapert2010,Khaneja2001,Zhang2011} and Magnetic Resonance Imaging~\cite{Conolly1986,Bernstein2004,Lapert2012}.
Recently, optimal control has attracted attention in view of applications to quantum information processing, for example as a tool to implement high-fidelity quantum gates in minimum time~\cite{cat,Garon2013,Boozer2012}. Generally, algorithms can also be designed to account for experimental imperfections or constraints
related to a specific material or device \cite{cat}. The possibility of including such constraints renders optimal control theory more useful in view of
experimental applications and helps bridge the gap between control theory and control experiments.

The standard numerical optimal control algorithms based on an iterative procedure compute the control fields through many time propagations of the state of the system, which can be prohibitive for systems of large dimensions in terms of computational time. This numerical limit can be bypassed by making use of parallel computing~\cite{decompspace,horton,maday2007}. In the case the computational time is divided by the number of computers, the method is said to be fully efficient. This full efficiency can be viewed as the physical limit in terms of performance of a parallel algorithm. While in applied mathematics different techniques have been developed using space or time decomposition~\cite{decompspace,maday2007}, very little has been done in quantum mechanics. The exponential growth of the Hilbert space dimension with the system size makes this question even more crucial in order to simulate the dynamics of complex quantum systems. Note that quantum control computations can also be speeded up by the parallelization of matrix exponential algorithms \cite{parallmatrix1,parallmatrix2} and by parallelizing density operator time evolutions using minimal sets of pure states \cite{skinner}.




This paper is not aimed at proposing a new optimization approach, but rather at describing and studying a general framework, namely the Intermediate State Method (ISM), introduced in \cite{maday2007}, which uses a time-parallelization to speed up the computation of optimal control fields. We investigate the efficiency of ISM on three benchmark quantum control problems, ranging from the control of coupled spin systems and the control of molecular orientation to the control of Bose-Einstein condensates.
As a by-product, we show under which conditions ISM can be made fully efficient.

The paper is organized as follows. Section~\ref{sec:2} is dedicated to the description of the time-parallelization method. The numerical schemes involved in this algorithm are defined in Sec.~\ref{sec:3}. Numerical results on the control of spin systems, molecular orientation and Bose-Einstein condensates are presented in Sec.~\ref{sec:4}. Conclusion and prospective views are given in Sec.~\ref{sec:5}.



\section{The time-parallelization Method}\label{sec:2}

We first introduce the optimal control problem and we derive the corresponding optimality conditions. We consider pure quantum states and we assume that the time
evolution is coherent. Note that the formalism can be easily extended to mixed states or to the control of evolution operators~\cite{palao}. The control process is aimed at maximizing the transfer of population onto a target state, but modification of the algorithms
in view of optimizing the expectation value of an observable is straightforward. The dynamics of the quantum system is governed by the Hamiltonian $H$. The initial and target states are denoted by $|\psi_i\rangle$ and $|\psi_f\rangle$, respectively and the general state of the system at time $t$, by $|\psi(t)\rangle$.
The dynamics of the quantum system is governed by the Schr\"odinger equation:
\begin{equation}\label{eq:dir}
i \partial_t |\psi(t)\rangle = H(u(t))|\psi(t)\rangle,
\end{equation}
where $u(t)$ is the field to be determined. The control time $T$ is fixed. The objective of the control problem is to maximize the figure of merit $\J$
\begin{equation*}
\J[u]=\Re [\langle\psi(T)|\psi_{f}\rangle]-\frac\alpha2\int_0^T u(t)^2dt,
\end{equation*}
$\alpha$ being a positive parameter which expresses the relative weight between the projection onto the target state and the energy of the control field.
A necessary condition to ensure the optimality of $u$ is given by the cancellation of the gradient of $\J$ with respect to $u$ \cite{Khaneja2005}:
\begin{equation}\label{eq:gradJ}
\nabla\J[u](t)=-\alpha u(t)+\Im [\langle \chi(t) |\partial_{u(t)} H (u(t))|\psi(t)\rangle]=0,
\end{equation}
where $|\chi(t)\rangle $ is the adjoint state that satisfies
\begin{equation}\label{eq:adj}
i \partial_t |\chi(t)\rangle = H(u(t))|\chi(t)\rangle,
\end{equation}
with the final condition $|\chi(t=T)\rangle=|\psi_f\rangle$.

\begin{figure}
\begin{tikzpicture}[scale=0.8]
             \coordinate (init) at (0,0);
             \coordinate (cible) at (5,3);
			\coordinate (P1) at (1,{1/25*1^2+.2*sin(40*pi*1)+1/5*(3-1/25*5^2-.2*sin(40*pi*5))});
			\coordinate (P2) at (2,{1/25*2^2+.2*sin(40*pi*2)+2/5*(3-1/25*5^2-.2*sin(40*pi*5))});
			\coordinate (P3) at (3,{1/25*3^2+.2*sin(40*pi*3)+3/5*(3-1/25*5^2-.2*sin(40*pi*5))});
			\coordinate (P4) at (4,{1/25*4^2+.2*sin(40*pi*4)+4/5*(3-1/25*5^2-.2*sin(40*pi*5))});
\begin{scope}[yshift=2.32cm,xshift=2.5cm,scale=.15,rotate=200]
\draw [fill] (0,1) -- (2,0) -- (0,-1)-- (.5,0) -- (0,1);
\end{scope}
\begin{scope}[yshift=0.10cm,xshift=2.5cm,scale=.15,rotate=33]
\draw [fill] (0,1) -- (2,0) -- (0,-1)--(.5,0) -- (0,1);
\end{scope}
\begin{scope}[yshift=-.85cm,xshift=7.5cm,scale=.1,rotate=230]
\draw [fill] (0,1) -- (2,0) -- (0,-1)-- (.5,0) -- (0,1);
\end{scope}
\begin{scope}[yshift=-2.04cm,xshift=7.5cm,scale=.1,rotate=52]
\draw [fill] (0,1) -- (2,0) -- (0,-1)-- (.5,0) -- (0,1);
\end{scope}
              \draw (init)  node {\huge $\bullet$} (cible)  node {\huge $\bullet$}  ;
         \draw[->,>=latex] (0,-.5) -- (5.5,-.5) ;
         \draw (0,-.6) -- (0,-.4) (1,-.6) -- (1,-.4) (2,-.6) -- (2,-.4) (3,-.6) -- (3,-.4) (4,-.6) -- (4,-.4) (5,-.6) -- (5,-.4)  ;

\draw(-.2,-.8)  node {\footnotesize $t=0$} ;
\draw(5.2,-.8)  node {\footnotesize $t=T$} ;

 \draw[->,>=latex] (3+.2,1.4) to[bend left=30] (7,{1/25*3^2+.2*sin(40*pi*3)+3/5*(3-1/25*5^2-.2*sin(40*pi*5))+.2-2}) ;
 \draw[->,>=latex] (6.1,-3) to[bend left=40] (2.5,-2.2) ;

         \draw[myred] (7.2,1.5)  node {\small $\left(|\varphi_n^u\rangle\right)_{n=1,\cdots,5}$} ;
         \draw (-.2,.6)  node {\small \it $|\psi_i\rangle$} ;
         \draw (5.85,3)  node {\small \it $|\psi_f\rangle$} ;
         \draw (4,0.2)  node {\small \it $|\psi (t)\rangle$} ;
         \draw (2,2.7)  node {\small \it $|\chi(t)\rangle$} ;

\draw plot [domain=0:5] (\x , {1/25*\x^2+.2*sin(40*pi*\x) });
\draw plot [domain=0:5] (\x , {1/25*\x^2+.2*sin(40*pi*\x) +3-1/25*5^2-.2*sin(40*pi*5)});
\draw[densely dashed] plot [domain=0:5] (\x , {1/25*\x^2+.2*sin(40*pi*\x)+\x/5*(3-1/25*5^2-.2*sin(40*pi*5))});

\draw[myred](P1) node {\Large $\bullet$} ;
\draw[myred](P2) node {\Large $\bullet$} ;
\draw[myred](P3)node {\Large $\bullet$} ;
\draw[myred](P4) node {\Large $\bullet$} ;

\draw (2-.2,{1/25*2^2+.2*sin(40*pi*2)+2/5*(3-1/25*5^2-.2*sin(40*pi*5))-.2}) rectangle
             (3+.2,{1/25*3^2+.2*sin(40*pi*3)+3/5*(3-1/25*5^2-.2*sin(40*pi*5))+.2}) ;

\begin{scope}[yshift=-4cm,xshift=2.5cm,scale=2]
\draw (2-.2,{1/25*2^2+.2*sin(40*pi*2)+2/5*(3-1/25*5^2-.2*sin(40*pi*5))-.75}) rectangle
             (3+.2,{1/25*3^2+.2*sin(40*pi*3)+3/5*(3-1/25*5^2-.2*sin(40*pi*5))+.2}) ;
\begin{scope}
\clip (2,{1/25*2^2+.2*sin(40*pi*2)+2/5*(3-1/25*5^2-.2*sin(40*pi*5))-.1})  rectangle
           (3,{1/25*3^2+.2*sin(40*pi*3)+3/5*(3-1/25*5^2-.2*sin(40*pi*5))+.1});
\draw plot [domain=0:5,samples=200] (\x , {1/11*\x^2-.1*sin(160*pi*\x)
                                                                                   -  1/11*2^2+.1*sin(160*pi*2)
                                                                          + 2/5*(3-1/25*5^2-.2*sin(40*pi*5))  });
\draw plot [domain=0:5,samples=200] (\x , {1/11*\x^2-.1*sin(160*pi*\x)
                                                                                   -  1/11*3^2+.1*sin(160*pi*3)
                                                                          + 1/25*3^2+.2*sin(40*pi*3)+3/5*(3-1/25*5^2-.2*sin(40*pi*5))  });

\end{scope}
\draw[myred](2,{1/25*2^2+.2*sin(40*pi*2)+2/5*(3-1/25*5^2-.2*sin(40*pi*5))}) node {\Large $\bullet$}
                            (3,{1/25*3^2+.2*sin(40*pi*3)+3/5*(3-1/25*5^2-.2*sin(40*pi*5))}) node {\Large $\bullet$} ;
\draw[->,>=latex] (1.9,.6) -- (3.15,.6) ;
\draw (2,.7) -- (2,.5) (3,.7) -- (3,.5) ;
\draw(2,.3)  node {\small $t_n$} ;
\draw(3,.3)  node {\small $t_{n+1}$} ;
\draw(2.9,.9)  node {\small $|\psi_n(t)\rangle $} ;
\draw(2.2,1.7)  node {\small $|\chi_n(t)\rangle $} ;
\end{scope}


         \draw[->,>=latex] ( .5,-.6) -- ( .5,-1);
			\node[draw,text width=.5cm,text centered,myred] at(.5,-1.6){\tiny \it CPU 1};
         \draw[->,>=latex] (1.5,-.6) -- (1.5,-1);
			\node[draw,text width=.5cm,text centered,myred] at(1.5,-1.6){\tiny \it CPU 2};
         \draw[->,>=latex] (2.5,-.6) -- (2.5,-1);
			\node[draw,text width=.5cm,text centered,myred] at(2.5,-1.6){\tiny \it CPU 3};
         \draw[->,>=latex] (3.5,-.6) -- (3.5,-1);
			\node[draw,text width=.5cm,text centered,myred] at(3.5,-1.6){\tiny \it CPU 4};
         \draw[->,>=latex] (4.5,-.6) -- (4.5,-1);
			\node[draw,text width=.5cm,text centered,myred] at(4.5,-1.6){\tiny \it CPU 5};

\end{tikzpicture}
\caption{Schematic description of the Intermediate state method (see the text for details).}\label{fig0}
\end{figure}
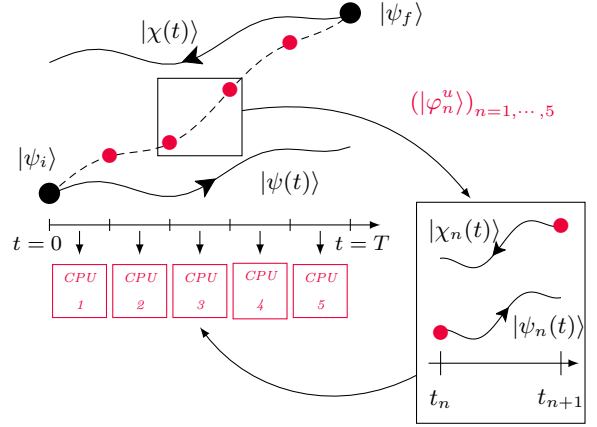

We now present ISM. A schematic description is displayed in Fig.~\ref{fig0}.
The main idea consists in considering a combination of the trajectories followed by $|\psi(t)\rangle$ and $|\chi(t)\rangle$~\cite{ST}. Given $N \in\N$, we decompose the interval $[0,T]$ into a partition of sub-intervals $[0,T]=\cup_{n=0}^{N-1} [t_n,t_{n+1}]$, with
$0=t_0<\cdots<t_N=T$. The parallelization strategy is based on this decomposition.
We consider an arbitrary control $u$ and we introduce
the sequence $|\phi^{u}\rangle=(|\phi^{u}_{n}\rangle)_{n=0,\cdots,N}$ that interpolates the state and adjoint state trajectories at time $t_n$ as follows:
\begin{equation}\label{eq:cible}
|\phi^{u}_{n}\rangle = \frac{T-t_n}{T}|\psi(t_n)\rangle+\frac{t_n}{T}|\chi(t_n)\rangle,
\end{equation}
where $|\psi(t)\rangle $ and $|\chi(t)\rangle $ are defined by Eqs.~\eqref{eq:dir} and \eqref{eq:adj}, respectively.
Note that $|\phi^{u}\rangle$ does not sample any usual dynamics, e.g. $|\phi^{u}\rangle$ does not correspond to a solution of Eq.~\eqref{eq:dir}.
Its initial and final states are $|\phi^{u}(0)\rangle =|\psi_0\rangle$ and $|\phi^{u}(T)\rangle =|\psi_f\rangle$, respectively. The choice of intermediate states made in Eq.~\eqref{eq:cible} is crucial to demonstrate Theorem \ref{th:1} below \cite{maday2007}.

We then introduce in each sub-interval the optimal control problem $\max_{u_n}\J_n[|\phi^{u}\rangle,u_n]$ defined by the maximization of the sub-functional:
\begin{eqnarray*}
& & \J_n[u_n,|\phi^{u}\rangle]= -\frac{1}{2} || |\psi_n(t_{n+1})\rangle -|\phi_{n+1}^u\rangle ||^2 \\
& & -\frac{\alpha_n}2\int_{t_n}^{t_{n+1}} u_n(t)^2dt,
\end{eqnarray*}
with $0\leq n\leq N-1$. In this problem, the state $|\psi_n\rangle$ is defined on $[t_n,t_{n+1}]$ by:
\begin{equation}\label{eq:dirpara}
i \partial_t |\psi_n(t)\rangle = H(u(t))|\psi_n(t)\rangle,
\end{equation}
starting from $|\psi_n(t_n)\rangle=|\phi^{u}_n\rangle$.
The penalization coefficient is defined by
$\alpha_n=\frac{t_{n+1}-t_n}{T}\alpha$.
Since
\begin{eqnarray*}
& & \J_n[u_n,|\phi^{u}\rangle]=-\frac{1}{2}|| |\psi_n(t_{n+1})\rangle ||^2-\frac{1}{2}|| |\phi_{n+1}^u\rangle ||^2 \\
& & +\Re[\langle \psi_n(t_{n+1})|\phi^{u}_{n+1}\rangle]-\frac{\alpha_n}2\int_{t_n}^{t_{n+1}} u_n(t)^2dt,
\end{eqnarray*}
maximizing $\J_n$ with respect to $u$ is equivalent to maximize a figure of merit of the form:
$\Re[\langle \psi_n(t_{n+1})|\phi^{u}_{n+1}\rangle]-\frac{\alpha_n}2\int_{t_n}^{t_{n+1}}
 u_n(t)^2dt$. In this way, each sub-problem has the same structure as the initial one.


We now review some properties of the time decomposition in order to establish the relation with the original optimal control problem.
Given an arbitrary trajectory $|\phi(t)\rangle$, we define an auxiliary figure of merit:
$$
\J_{\parallel}[u,|\phi\rangle]=\sum_{n=0}^{N-1}\beta_{n}\J_{n}[u_n,|\phi\rangle ],
$$
with $\beta_n=\frac{T}{t_{n+1}-t_n}$. A first relation between $\J_{\parallel}$ and $\J$ is given in Theorem \ref{th:1} (see Ref.~\cite{maday2007} for the proof).
\begin{theorem}\label{th:1}
Given an arbitrary control $u$, we have:
$$
|\phi^{u}\rangle={\rm argmax}_{|\phi\rangle}\big(\J_{\parallel}[u,|\phi\rangle]\big).
$$
Moreover, the following relation is satisfied:
$$
\J_{\parallel}[u,|\phi^{u}\rangle]=\J[u].
$$
\end{theorem}
As a by-product, this Theorem allows us to compute in parallel $\J[u]$, knowing only the sequence  $|\phi_u\rangle$.
A similar relation also holds between the gradients of the functionals, as stated in Theorem \ref{th:gradient}.
\begin{theorem}\label{th:gradient}
Given an arbitrary control $u$, we have:
$$\nabla \J [{u}]_{|[t_{n},t_{n+1}]} =\beta_{n}\nabla \J_n [{u_{|[t_{n},t_{n+1}]}},|\phi^{u}\rangle].$$
\end{theorem}
This result provides a new interpretation of the time-parallelized method since the sequence $|\phi^u(t_n)\rangle$, $n=0,\cdots,N$ of intermediate states enables the decomposition of the computation of the gradient.\\

\noindent {\it Proof:} Let us consider a fixed value $n$, with $0\leq n\leq N-1$, $t\in [t_n,t_{n+1}]$ and denote by $|\chi_n(t)\rangle $ and $|\psi_n(t)\rangle$ the trajectories defined by
\begin{equation*}
i \partial_t |\chi_n(t)\rangle = H(u(t))|\chi_n(t)\rangle,
\end{equation*}
and
\begin{equation*}
i \partial_t |\psi_n(t)\rangle = H(u(t))|\psi_n(t)\rangle,
\end{equation*}
with $|\chi_n(t_{n+1})\rangle=|\phi^u_{n+1}\rangle $ and $|\psi_n(t_{n})\rangle=|\phi^u_{n}\rangle $. For $t\in [t_n,t_{n+1}]$, we repeat with $\J_n$ the computation made to derive the gradient of $\J$:
\begin{eqnarray*}
\nabla \J_n [u_{|[t_{n},t_{n+1}]}, |\phi^{u}\rangle ](t) &=& \Im\left(\langle \chi_n(t)|\partial_{u(t)}H|\psi_n(t)\rangle \right)\\
&&-\alpha_n u(t).
\end{eqnarray*}
Using the fact that:
\begin{equation*}
|\chi_n(t)\rangle=\frac{(T-t_{n+1})}{T}|\psi(t)\rangle + \frac{t_{n+1}}{T} |\chi(t)\rangle,
\end{equation*}
and
\begin{equation*}
|\psi_n(t)\rangle=\frac{(T-t_{n})}{T}|\psi(t)\rangle + \frac{t_{n}}{T} |\chi(t)\rangle,
\end{equation*}
we arrive at:
\begin{eqnarray*}
\nabla \J_n [u_{|[t_{n},t_{n+1}]},|\phi^{u}\rangle](t) &=& \frac{1}{\beta_n}\Im\left(\langle \chi(t)|\partial_{u(t)}H |\psi(t)\rangle \right)\\
\\&&-\frac{\alpha}{\beta_n}u(t),
\end{eqnarray*}
and the result follows.\qed \\

We now give the general structure of ISM. Let $1\geq \eta >0$ and $u^{(0)}$ be an initial control field.
\begin{algorithm}~\label{alg:ISM}
\begin{enumerate}
\item Set $Err=1$, $k=0$.
\item While $Err > \eta$, do:\\
\begin{enumerate}
\item\label{etape:1} Compute on $[0,T]$ the trajectories $|\psi^{(k)}(t)\rangle$ and $|\chi^{(k)}(t)\rangle$ associated with $u^{(k)}$ according to Eqs.~\eqref{eq:dir} and~\eqref{eq:adj}.
\item\label{etape:2} Compute $|\phi^{(k)}(t)\rangle=|\phi^{u^{(k)}}(t)\rangle$ according to Eq.~\eqref{eq:cible}.
\item\label{etape:3} On each sub-interval $[t_n,t_{n+1}]$ compute in parallel an approximation of the solution $u^{(k+1)}_n$ of the problem
  $\max_{u_n} \J_n[u_n,|\phi_k\rangle]$.
\item\label{etape:4} Define $u^{(k+1)}$ as the concatenation of the controls $u^{(k+1)}_n$, $n=1,\cdots,N-1$.
\item Set $Err=\sum_{n=0}^{N-1}\int_{T_{n}}^{T_{n+1}}\|\nabla \J_n [u^{(k+1)}_{|[t_{n},t_{n+1}]}](t)\|dt$.
\item Set $k=k+1$.
\end{enumerate}
\end{enumerate}
\end{algorithm}
Step~\ref{etape:1} contradicts the parallelization paradigm, since it requires a sequential solving of an evolution equation on the full interval $[0,T]$. We will see how this problem can be solved in Sec.~\ref{sec:4}. However, note that the most time consuming step, namely Step~\ref{etape:3}, is achieved in parallel.
\section{Description of the numerical methods used in the parallelization}\label{sec:3}

Different schemes can be used to implement the time-parallelized algorithm in practice. This requires two ingredients, a numerical scheme to solve approximately the evolution equations of Steps~\ref{etape:1} and~\ref{etape:3} and an optimization procedure for the sub-problem of Step~\ref{etape:3}.
In this paragraph, we give some details about the used numerical methods and we explain how the full efficiency can be approached in the case of quantum systems of sufficiently small dimensions.

In the different numerical examples, we consider two numerical solvers for the Schr\"odinger Equation~\eqref{eq:dir}: a Crank-Nicholson scheme and a second order Strang operator splitting.
Such solvers can be described through an equidistant time-discretization grid $t_a=t_0<t_1<\cdots<t_J=t_b$ of an interval $[t_a,t_b]$.
The time step is denoted by $\tau = (t_b-t_a)/(J-1)$ for some $J\in\N$. For each time grid point $t_j$, we introduce the state $|{\psi}_j\rangle $ and the control $u_j$, which are some  approximations of the exact state $|\psi(t_j)\rangle $ and of the exact control field $u(\frac{t_j+t_{j+1}}{2})$.

The Crank-Nicholson algorithm is based on the following recursive relation:
\begin{equation}\label{eq:CN}
  \frac{i}{\tau}(|{\psi}_{j+1}\rangle-|{\psi}_{j}\rangle)= \frac{H(u_j)}{2}(|{\psi}_{j+1}\rangle+|{\psi}_{j}\rangle),
\end{equation}
which can be rewritten in a more compact form as
\begin{equation}
  (\openone+L_j)|\psi_{j+1}\rangle  = \big(\openone-L_j\big)|\psi_j\rangle\,,
  \label{eq:CN-Ln}
\end{equation}
where $\openone$ is the identity operator and $L_n := i\frac{\tau}{2}H(u_{n}(t))$.

The second order Strang operator splitting is rather used in the case of infinite dimensional systems. Indeed, this method
is particularly relevant when the Hamiltonian includes a differential operator. We consider, for example, the case $H(u(t))=-\Delta + V(u(t),x)$ where
$\Delta$ denotes the Laplace operator and $V(u(t))=V(u(t),x)$ is a scalar potential. In this case, Strang's method gives rise to the iteration:
\begin{equation}\label{eq:OS}
  |{\psi}_{j+1}\rangle= \exp(-\frac{i\tau}2\Delta)\exp(-i\tau V(u_j))\exp(-\frac{i\tau}2\Delta)|{\psi}_{j}\rangle.
\end{equation}
In Eq.~\eqref{eq:OS}, each product can be determined very quickly since the operator $V(u_j)$ is diagonal in the physical space, while $\Delta$ is
generally diagonal in the Fourier space, and the change of basis can be achieved efficiently by fast Fourier transform.

These schemes provide a second order approximation
with respect to time, which leads to an accurate approximation of the trajectory $|\psi(t)\rangle$.
In addition, both propagators automatically preserve the normalization of the wave function, which is very interesting to avoid non-physical solutions.
A specific advantage of these solvers is that they allow an exact differentiation with respect to the control in the discrete setting, in the case of scalar control for
the Strang solver~\eqref{eq:OS} and in any case for the Crank-Nicholson solver~\eqref{eq:CN}.

We now explain how the full efficiency can be reached with the parallelization algorithm. Both solvers lead to a linear relation between the initial and the final states of the system of the form:
\begin{equation*}
|{\psi}_{J}\rangle = M(u) |{\psi}_{0}\rangle.
\end{equation*}
As an example, for the Crank-Nicholson solver, we have:
\begin{equation*}
~M(u)=\Pi_{j=0}^{J-1}(\openone+L_j)^{-1}(\openone-L_j).
\end{equation*}
The matrix $M(u)$ can be computed in parallel during the propagation of Eq.~\eqref{eq:CN-Ln}. Knowing the state $|{\psi}_{0}\rangle$,
this matrix enables to
compute in one matrix-vector product $|{\psi}_{J}\rangle$. As a consequence, this propagator assembling technique
allows us to avoid the sequential solving on the full interval $[0,T]$ in Step~\ref{etape:1} of Algorithm~\ref{alg:ISM}.
More precisely, assume for example that at iteration $k$ of Algorithm~\ref{alg:ISM} and on each sub-interval, a
matrix $M_n(u^{(k)}_n)$ is computed and transmitted to the main processor. The computations
of the sequences $|\psi(t_n)\rangle$ and $|\chi(t_n)\rangle$, that are required to define the intermediate states
$|\phi^{(k)}(t_n)\rangle$, can be achieved in $2N$ matrix-vector products. Due to storage and communications issues of the matrices, note that this approach
can only be used for quantum systems of small dimensions.

We conclude this paragraph by presenting a way to derive the gradient of time-discretized figures of merit of the form:
$$\J_{\tau}[u]= \Re\langle {\psi}_J|\psi_f\rangle-\frac{\alpha}{2} \tau \sum_{j=0}^{J-1} u_j^2.$$
We consider the case of a Crank-Nicholson solver, but similar computations can be made for Strang's solver.
We introduce the functional $\L_{\tau}$ defined by:
\begin{eqnarray*}
& \L_{\tau}[u,|\psi\rangle,|\chi\rangle]=\J_{\tau}[u]+ \Re (\sum_{j=0}^{J-1} \langle \chi_j|\openone+L_j|\psi_{j+1}\rangle  \\
& - \langle \chi_j|\openone-L_j|\psi_j\rangle).
\end{eqnarray*}
Since $L_j$ is anti-hermitian, differentiating $\L_{\tau}$ with respect to $|\psi\rangle$ gives rise to the discrete adjoint evolution equation:
\begin{equation*}
(\openone-L_{j-1})|\chi_{j-1}\rangle  = (\openone+L_j)|\chi_j\rangle,
\end{equation*}
with the final condition $(\openone-L_{J-1})|\chi_{J-1}\rangle =|\psi_J\rangle $. To derive the gradient of $\J_{\tau}[u]$, it remains to
differentiate $L_{\tau}$ with respect to $u$, which leads to the $j$-th entry of the gradient of $\J_{\tau}$:
\begin{equation*}
\left(\nabla \J_{\tau}[u]\right)_j= \alpha dt u_j + \frac{i \tau}2 \langle \chi_j|\partial_u H|\psi_{j+1}+\psi_j\rangle.
\end{equation*}

In the sequel, we use this result to implement a constant step gradient method: the approximation of the solution of the sub-problem in Step~\ref{etape:4} is computed by
iterating on $\ell$ in the formula:
\begin{equation}\label{gradStep}
 u^{\ell+1}=u^\ell - \rho\nabla \J_{\tau}(u^\ell),
\end{equation}
for some $\rho>0$.\\

Other optimization methods such as pseudo or quasi-Newton approaches can be used to perform Step~\ref{etape:3}.


\section{Numerical results}\label{sec:4}
This section is dedicated to some numerical results obtained with ISM, used with the schemes presented in Sec.~\ref{sec:3}. The efficiency of this approach is illustrated on three benchmark examples in quantum control \cite{Brif2010,cat}, namely the control of a system of coupled spins, the control of molecular orientation and the control of a Bose-Einstein condensate whose dynamics is governed by the Gross-Pitaevskii equation.
\subsection{Control of a system of five coupled spins}\label{sec:5:1}
In this paragraph, we consider the control of a system of coupled spin 1/2 particles. The principles of control in Nuclear Magnetic Resonance being described in different books \cite{Bernstein2004,Ernst1987,Levitt2008}, only a brief account will be given here in order to introduce the used model. We investigate the control of a system of coupled spins by means of different magnetic fields acting as local controls on each spin. Each field only acts on one spin and does not interact with the others, i.e. the spins are assumed to be selectively addressable.

We introduce a system of 5 coupled spins \cite{spin5_marx1,spin5_marx2}, the evolution of which is described by the following Hamiltonian:
\begin{equation}
H=H_0+\sum_{k=1}^5[u_x^{k}(t)I_x^{(k)}+u_y^{k}(t)I_y^{(k)}]
\end{equation}
where the operators $I_x^{(k)}$ and $I_y^{(k)}$ are, up to a factor, Pauli matrices which only act on the $k$th spin:
\begin{equation*}
I_x:=\left(\begin{array}{cc}0&1/2\\1/2&0\end{array}\right),  I_y:=\left(\begin{array}{cc}0&-i/2\\i/2&0\end{array}\right),
\end{equation*}
\begin{equation*}
I_z:=\left(\begin{array}{cc}1/2&0\\0&-1/2\end{array}\right).
\end{equation*}
We assume that the free evolution Hamiltonian $H_0$ is associated with a topology \cite{spin5_marx1,spin5_marx2} defined by:
\begin{eqnarray*}
& & H_0	= 2\pi(J_{12} I_z^{(1)}I_z^{(2)}+J_{13}I_z^{(1)}I_z^{(3)}+\\
& &+J_{23} I_z^{(2)}I_z^{(3)}+J_{25}I_z^{(2)}I_z^{(5)}+J_{34} I_z^{(3)}I_z^{(4)}).
\end{eqnarray*}
Note that this model system is valid in heteronuclear spin systems if the coupling strength between the spins is small with respect to the frequency shifts \cite{Ernst1987,Levitt2008}. The coupling constant between the spins is taken to be uniform and equal to $J_p=140$. For the numerical simulations, we move to the density matrix formalism with $I_x^{(1)}$ and $I_x^{(5)}$ as initial and final states, respectively. The control time is fixed to $T=J_p/10$. The parameter $\alpha$ is set to 0.

\begin{figure}
\includegraphics[width=.75\linewidth]{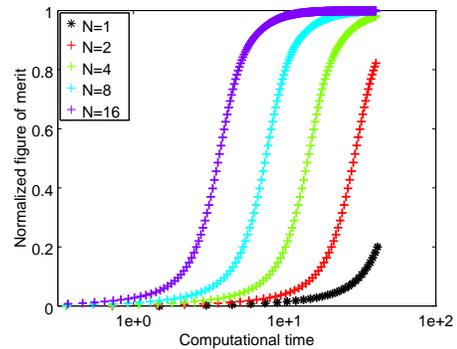}
\caption{Evolution of the normalized figure of merit at each iteration and for different values of $N$ (the number of processors) with respect to computational time (wall-clock time) in the case of the control of a spin system. As stated in Theorem~\ref{th:gradient},
the values obtained after a given number of iterations are the same for all values of $N$. Log-scale is used in the x-axis.}
\label{fig1}
\end{figure}

For the time-parallelization, we consider a uniform grid, so that $t_{n+1}-t_n=T/N$, for $n=0,\cdots,N-1$ and we compare the results for different values of $N$. The time discretization is done by the Crank-Nicholson method of Eq.~\eqref{eq:CN} with the time step $T/2^{15}$. In Step~\ref{etape:3}, one iteration of the constant step gradient descent method [see Eq.~\eqref{gradStep}] is used, with $\rho=10^4$. As a result of Theorem~\ref{th:gradient},
the values obtained after a given number of iterations are the same for all values of $N$.
In this way, the method is almost fully efficient. It is actually equivalent
to a standard gradient method, except that the computation of the gradient is done in parallel. The computational effort is therefore exactly divided by the number of processors, and the full efficiency is only limited by the memory usage and also by the communication between processors required by the update of the intermediate states in Steps~\ref{etape:1} and~\ref{etape:2} of the algorithm.
Figure~\ref{fig1} displays the figure of merit with respect to the parallel computational time.

\begin{figure}[!htpb]
\centering
      \includegraphics[width=\linewidth]{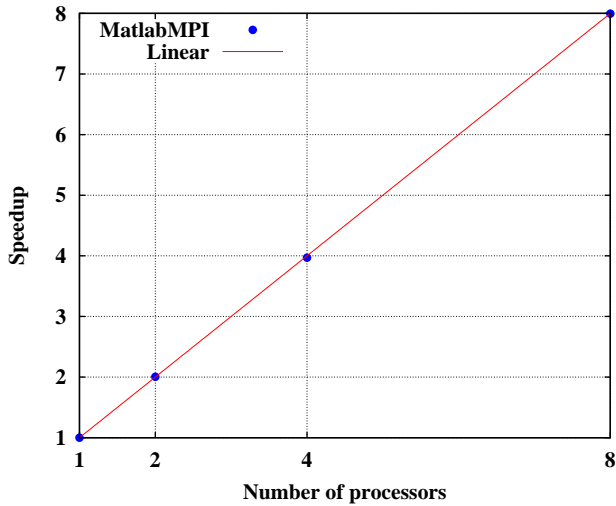}
\caption{Speedup $S(1.74,N)$ of the parallel implementation (y- axis) with respect to the number of processors $N$ (x- axis). The blue dots indicate the speedup achieved with MatlabMPI for $\varepsilon=1.74$ (see Table \ref{tabefficiency} for details). The red solid line corresponds to a linear evolution of the speedup as a function of $N$.}
\label{figspeedup}
\end{figure}

\begin{figure}[!htpb]
\centering
\begin{tabular}{cc}
      \includegraphics[scale =0.35]{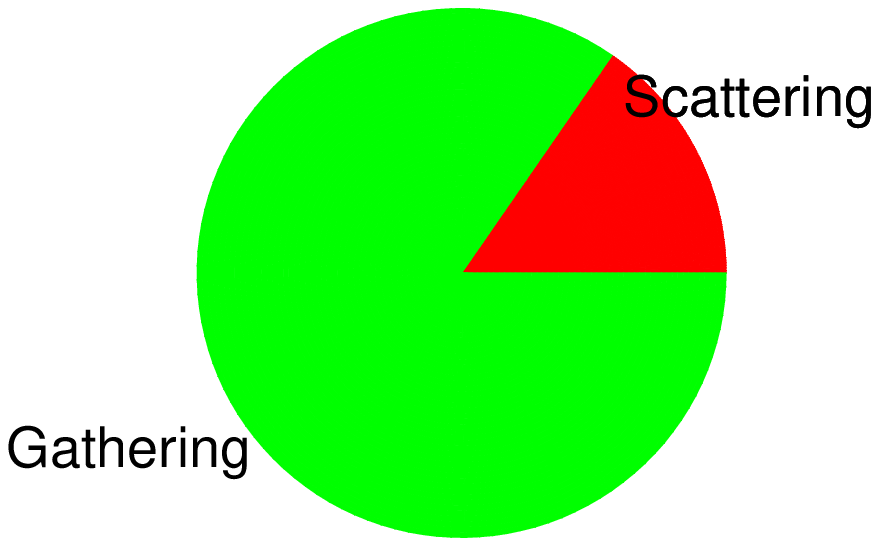} &
      \includegraphics[scale =0.35]{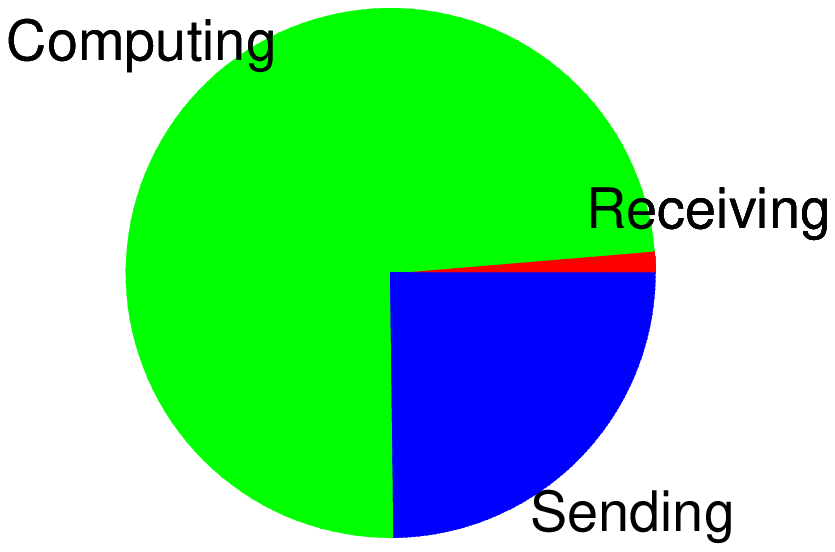}\\
 Master processor &       Slave processor\end{tabular}
 \caption{MPI profiling of the tasks performed by the master and a slave processors with $N=4$ and for 20 iterations of the optimization process. The proportions are computed separately regarding the own wall-clock timing of a given processor.}\label{profiling}
\end{figure}

In order to evaluate more precisely the efficiency of the algorithm, we give some details about the speedup of the numerical computations. Numerical simulations are implemented with Matlab, where the parallelization is realized using the open source library MatlabMPI~\cite{matlabmpi}. The  tests have been carried out on a shared memory machine under a linux system with a core of Intel(R) Xeon(R) CPU type (@ 2.90GHz with 198 Giga byte shared memory). The parallel computation uses $N$  processors where $N$ stands, as above, for the number of sub-intervals of the time domain decomposition. In Fig.~\ref{figspeedup}, parallel numerical performances are compared with the sequential performance which is obtained when a single processor is used to treat the whole time domain.
Given $\varepsilon>0$, we introduce the quantities $S(\epsilon,N):={t(\epsilon,1)}\slash{t(\epsilon,N)}$ and $Eff(\epsilon,N):=10^2({S(\epsilon,N)}\slash{N}),$  as the parallel speedup and the efficiency respectively, where $t(\epsilon,N)$ denotes the computational time (with $N$ processors) necessary to reach a value $\mathcal{J}[u^k]$ such that  $\mathcal{J}[u^\infty]-\mathcal{J}[u^k]<\varepsilon$, where $u^\infty$ is
the value of the sequence $u^k$ obtained at the numerical convergence. Figure~\ref{figspeedup} displays results about the speedup of the parallel implementation.
We observe that the algorithm behaves as expected when increasing the number of processors. Despite the use of Input/Output (I/O) data files to ensure the communication between CPUs (as required by MatlabMPI), ISM achieves a linear scalability. A profiling of the parallel computing is reported in Fig.~\ref{profiling} where we present the time spent to achieve the communications for the master and one slave processors during 20 iterations of the optimization process. As expected, the speedup is independent on the value of $\epsilon$.  This point is clearly exhibited in Table~\ref{tabefficiency}. The communications in MatlabMPI are of point-to-point type through I/O files, a for-loop is therefore necessary to cover all sender and receiver processors. In this view, a slave processor waits for its turn in order to be able to read the message from the master processor. On the contrary, the printing message addressed to the master processor is done on slave processors and hence is non-blocking. 

Figure~\ref{profiling} shows that, in the case $N=4$, the slave processor mainly works on parallel computing. Its communication part is shared between sending and receiving data. The sending part is longer due to the amount of data to treat, which consists not only of partial control but also of propagator matrices.

\begin{table}[htbp]
\begin{tabular}{|c|c|c|c|c|c|}\hline
$N$ processors & 1 & 2 & 4 & 8 \\\hline
$Eff(0.32,N)$ & 100\% & 100.5\%&99.7 \%& 100.9\%\\\hline
$Eff(0.89,N)$ & 100\%&100.2 \%&99.2 \%& 99.9\%\\\hline
$Eff(1.74,N)$ & 100\% & 100.2\%& 99.2\%& 99.9\%\\\hline
\end{tabular}
\caption{Efficiency $Eff(\epsilon,N)$ of the parallel MPI implementation. Note the full efficiency performance through several snapshots in the running parallel computing.}\label{tabefficiency}
\end{table}

The full efficiency of ISM is clearly shown in the numerical simulations, see Table.~\ref{tabefficiency}. 
Note that, in some cases, the efficiency is greater than $100\%$ because the full problem requires more memory, and thus spends more time in hardware storage processes. On the contrary, the parallel computing uses a smaller amount of data. Also, depending on the processor architecture, the computational time shall behave nonlinearly with respect to the size of the data.
Similar linear and super linear speedup have been observed in Ref.~\cite{kepner2004matlabmpi} with MatlabMPI.

\subsection{Optimal control of molecular orientation}
In a second series of numerical tests, we consider the control of molecular orientation by THz laser fields. Molecular orientation \cite{revieworientation1,revieworientation2} is nowadays a well-established topic both from the experimental \cite{orientation2color,orientationTHZ} and theoretical points of views \cite{dionorientation,averbukhorientation,daemsorientation,lapertorientation,tehini2color}. Different optimal control analyses have been made on this quantum system \cite{DST,lapertmono,ohtsukiorientation1,ohtsukiorientation2}.

In this paragraph, we consider the control of a linear polar molecule, HCN, by a linearly polarized THz laser field $E(t)$. We assume that the molecule is in its ground vibronic state and described by a rigid rotor. In this case, the Hamiltonian of the system can be written as:
$$
H(t)=BJ^2-\mu_0\cos\theta E(t)-\frac{E(t)^2}{2}[(\alpha_\parallel-\alpha_\perp)\cos^2\theta+\alpha_\perp],
$$
where $B$ is the rotational constant, $J^2$ the angular momentum operator, $\mu_0$ the permanent dipolar moment, $\alpha_\parallel$ and $\alpha_\perp$ the dipole polarizability components parallel and perpendicular to the molecular axis, respectively, and $\theta$ the angle between the molecular axis and the polarization direction of the electric field. At zero temperature, the dynamics of the system is ruled by the following differential equation:
$$
i\frac{\partial|\psi(t)\rangle }{\partial t}=H(t)|\psi(t)\rangle
$$
where the initial state at $t=0$ is $|0,0\rangle$, in the basis of the spherical harmonics $\{|j,m\rangle,~j\geq 0,~-j\leq m\leq j\}$. Numerical values are taken to be $B=6.6376\times 10^{−6}$, $\mu_0=1.1413$, $\alpha_\parallel=20.055$ and $\alpha_\perp=8.638$, in a.u. We refer the reader to Ref.~\cite{DST} for details on the numerical implementation of this control problem. The $\alpha$ parameter is chosen as a time-dependent function of the form $10^5(\frac{t-T/2}{T/2})^6+10^4$ in order to design a control field which is experimentally relevant~\cite{DST}. The target state is the eigenvector of the observable $\cos\theta$ with the maximum eigenvalue in the subspace such that $j\leq 4$ \cite{daemsorientation,targetorientation}.

\begin{figure}
\includegraphics[width=\linewidth]{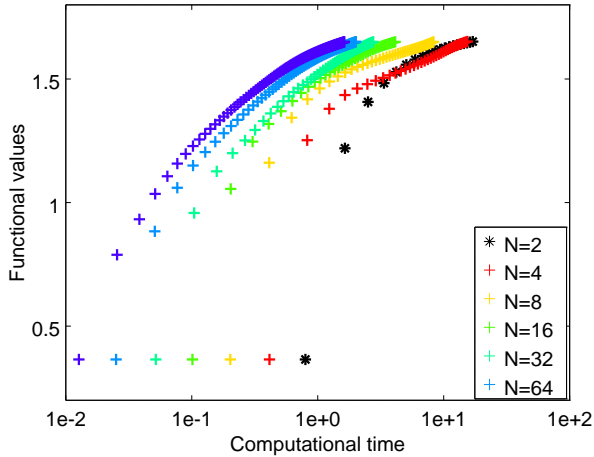}
\caption{Functional values at each iteration and for various values of $N$  with respect to computational time for the control of molecular orientation by a monotonic algorithm. Log-scale is used in the x-axis.}
\label{fig2}
\end{figure}

Having investigated the implementation issues in Sec.~\ref{sec:5:1}, we focus in this paragraph on the efficiency $Eff^\star(\epsilon,N)$ achieved when neglecting the time
associated with I/O communications. As a consequence, the results hereafter do not depend on the used computer and software. We study the efficiency of the parallelization method in the cases where the optimization solver consists in one step of either a monotonic algorithm or a Newton method. We start with a simulation using monotonic algorithm (see~\cite{S,Reich2012} for details about this method). Given a target value $\varepsilon$ of the figure of merit, we measure the computational time necessary to obtain it. In the numerical computation, we use $\epsilon=0.3$, while the optimum has been numerically estimated as $0.2909$~\cite{DST}. The values of the figure of merit are plotted in Fig.~\ref{fig2}.
In this case, the full efficiency is not obtained, as reported in Table~\ref{tabmono}.
\begin{table}
\begin{tabular}{|c|c|}
\hline $N$ & $Eff^\star(\varepsilon,N)$\\
\hline
\hline 1 &  100\% \\
\hline 2 &  55.2\% \\
\hline 4 &  51.2\%\\ 
\hline 8 &  51.5\%\\ 
\hline 16 & 38.5\%\\ 
\hline 32 & 26.4\%\\ 
\hline 64 & 16.3\%\\ 
\hline
\end{tabular}
\caption{Efficiency of ISM for various values of $N$ in the case of the control of molecular orientation with a monotonic algorithm.}
\label{tabmono}
\end{table}

In the different numerical simulations, we observe that this method seems to be significantly more efficient than gradient descent solvers.
The wall-clock computational time is a decreasing function of $N$, so that solving benefits from large parallelization.
However, $Eff^\star(\varepsilon,N)$ does not appear to be a monotonic function of $N$. The analysis of this point is out of the scope
of this paper. 

\begin{figure}
\includegraphics[width=\linewidth]{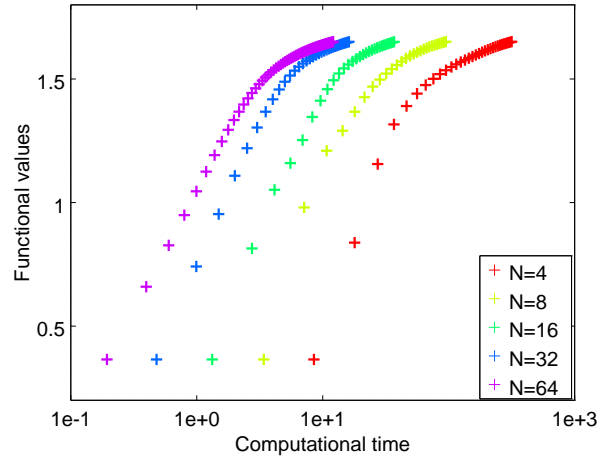}
\caption{Same as Fig~\ref{fig2} but for a Newton solver.}
\label{fig3}
\end{figure}

We repeat this test with one iteration of the Newton method as optimization solver.
More precisely, we implement a matrix-free version of the algorithm that
updates the control by means of a GMRES routine \cite{GMRES1,GMRES2}. In this case, we observe that this approach
actually enables to obtain convergence, the algorithm does not converge for $N=1$ and $N=2$ but converges for larger values of $N$.
The values of the figure of merit are plotted in Fig.~\ref{fig3}. Since the algorithm does not converge for $N=1$, $t(\varepsilon, 1)$ and $Eff^\star$ are not defined. In this case, we consider the quantity $ N\cdot t(\varepsilon, N)$ to measure the efficiency of the process. The results are presented in Table~\ref{tabnewt}.
\begin{table}
\begin{tabular}{|c|c|}
\hline $N$ & $N\cdot t(\varepsilon,N)$\\
\hline\hline 1 &  - \\
\hline 2 &  - \\
\hline 4 & 1264.763737 \\
\hline 8 & 759.976361 \\
\hline 16 & 589.424517 \\
\hline 32 & 516.603943 \\
\hline 64 & 774.557304 \\\hline
\end{tabular}
\caption{Equivalent sequential time of ISM in the case of a Newton solver for the control of molecular orientation.}
\label{tabnewt}
\end{table}
\subsection{Optimal control of Bose-Einstein condensates}~\label{sec:nonlin}

The last example investigated in this work deals with the optimal control of Bose-Einstein condensates \cite{BECfolman}. This subject has been extensively studied in the past few years \cite{Alfio,BECborzi,BECbucker,BECjager,BEClangen,BEClapert,BECschaff}. Following Ref.~\cite{Alfio}, we consider the control of a condensate in magnetic microtraps whose dynamics is ruled by the Gross-Pitaevskii equation:
\begin{equation}
\frac{\partial}{\partial t} |\psi (x,t)\rangle = \left(H_0+ V(x,\lambda(t))\right)|\psi (x,t)\rangle,
\end{equation}
where $|\psi(x,t)\rangle$ is the state of the system, $H_0=-\frac{1}{2}\frac{\partial^2}{\partial{x^2}}+\kappa \langle\psi (x,t)|\psi(x,t)\rangle$, $\lambda$ is the radio-frequency control field and $\kappa$ a positive coupling constant.
The potential $V$ is defined by
$$V(x,\lambda)=\left\{
\begin{array}{cl}
\frac 12 \left(|x|-\frac{\lambda d}{2} \right)^2 & {\rm for}\ |x|>\frac{\lambda d}{4} \\
\frac 12 \left(\frac{(\lambda d)^2}{8} - x^2\right)   & {\rm otherwise,}
\end{array}
\right. $$
with $d>0$. Unitless parameters are used here. We refer the reader to Ref.~\cite{Alfio} for details on the model system. The figure of merit associated with this control problem is:
\begin{equation}
{\cal J}[\lambda]= \Re[\langle\psi(T)|\psi_{f}\rangle ].
\end{equation}
The final and initial states of the control problem are respectively the ground states of the Hamiltonians $H_0+V(x,0)$ and $H_0+V(x,1)$. The parameter $\alpha$ is set to 0.

Due to the nonlinearity of the model system, the preceding approach has
to be adapted. In this way, we do not consider anymore the sequence of states
$|\varphi_n^u\rangle$ [see Eq.~\eqref{eq:cible}] but we
split up these intermediate states into two sets. The sequence of
initial states
is taken on the trajectory $|\psi(t)\rangle$, i.e. defined by
$|\psi(t_n)\rangle$, $n=0,\cdots,N$, while the sequence of target states  is taken on the adjoint trajectory
$|\chi(t)\rangle$, i.e. defined by $|\chi(t_n)\rangle$, $n=0,\cdots,N$.
The maximization problem in Step~\ref{etape:3} of
Algorithm~\ref{alg:ISM} is therefore replaced by the maximization of
the sub-functional:
\begin{equation}
\J_n[u_n,|\psi^{u}\rangle]= \Re[\langle
\psi_n(t_{n+1})|\chi(t_{n+1})\rangle ],
\end{equation}
with $0\leq n\leq N-1$. In this problem, the state $|\psi_n\rangle$ is
defined on $[t_n,t_{n+1}]$ by Eq.~\eqref{eq:dirpara},
but starting from $|\psi_n(t=t_n)\rangle=|\psi(t_n)\rangle$. The rest of
the procedure remains unchanged.


\begin{figure}
\includegraphics[width=\linewidth]{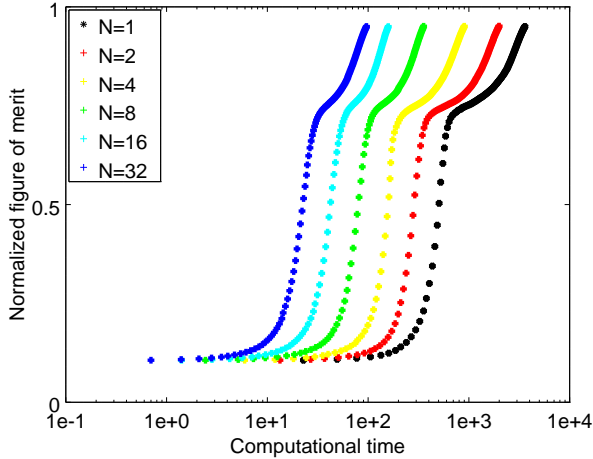}
\caption{Evolution of the normalized figure of merit at each iteration and for various values of $N$  with respect to computational time for the control of a Bose-Einstein condensate.}\label{fig4}
\end{figure}
This modification does not affect dramatically the computational time since these trajectories are not
computed sequentially but in parallel. We then use the propagator assembling technique presented in Sec.~\ref{sec:4} to compute the sequences of initial and final states. The numerical values of the parameters are set to $T=8$, $\kappa=1$ and $d=10$. The space domain we consider is $[-10,10]$.
For the space discretization, we consider a uniform grid composed of $50$ points. The time discretization is achieved
with a time grid of $2^9$ points, and we use Strang's splitting~\eqref{eq:OS} to compute the trajectories.
The optimization solver in Step~\ref{etape:3} consists in one iteration of the constant step gradient descent method, see Eq.~\eqref{gradStep}, with $\rho=10^{-1}$.

The results are presented in Fig.~\ref{fig4}.

\begin{table}
\begin{tabular}{|c|c|}
\hline $N$ & $Eff^\star(\varepsilon,N)$\\
\hline
\hline 1 &  100\%\\ 
\hline 2 &  90\%\\ 
\hline 4 &  99.7\%\\ 
\hline 8 &  126.3\%\\ 
\hline 16 & 141.8\%\\ 
\hline 32 & 116.65\%\\ 
\hline
\end{tabular}
\caption{Efficiency of ISM for various values of $N$ in the case of the control of a Bose-Einstein condensate.}
\label{tabGP}
\end{table}
We observe that the full efficiency is reached, as confirmed in Table~\ref{tabGP}. This feature is certainly a consequence of the nonlinear setting. The sub-control problems are simpler not only because of the size reduction induced by the time decomposition, but also because of the dynamics itself, which is simplified on a shorter time interval.


\section{Conclusion and perspectives}\label{sec:5}
In this work, we have investigated the numerical efficiency of a time parallelized optimal control algorithm on standard quantum control problems extending from the manipulation of spin systems and molecular orientation to the control of Bose-Einstein condensates. We have shown that the full efficiency can be reached in the case of a linear dynamics optimized by means of gradient methods. On the contrary, full efficiency is not achieved when using monotonic algorithms and Newton solvers. In the case of a Newton method,
the parallelization setting reduces the length of time intervals where the solver is used, and makes
the subproblems easier to solve. Such a property is also observed in the case of nonlinear dynamics, as shown with the example of Bose-Einstein condensates.
The results of this work can be viewed as an important step forward for the implementation of parallelization methods in quantum optimal control algorithms. Their use will become a prerequisite in a near future to simulate quantum systems of increasing complexity.

\noindent\textbf{ACKNOWLEDGMENT}\\
S.J. Glaser acknowledges support from the DFG (Gl 203/7-1), SFB 631 and the BMBF
FKZ 01EZ114 project. D. Sugny and S. J. Glaser acknowledge support from the ANR-DFG research program Explosys (ANR-14-CE35-0013-01; DFG-Gl 203/9-1). J.S was partially supported by the Agence Nationale
de la Recherche (ANR), Projet Blanc EMAQS number ANR-2011-BS01-017-01. This work has been done with the support of the Technische Universit\"at M\"unchen – Institute for Advanced
Study, funded by the German Excellence Initiative and the European Union Seventh Framework Programme under grant agreement 291763.


\end{document}